\documentclass[12pt]{article}
\usepackage{amsmath,amsfonts,amssymb}

\newcommand{\C}{\mathbb C}
\newcommand{\Q}{\mathbb Q}
\newcommand{\Z}{\mathbb Z}
\newcommand{\cB}{\mathcal B}
\newcommand{\PP}{\mathbb P}
\newcommand{\Oh}{\mathcal O}
\newcommand{\eps}{\varepsilon}
\DeclareMathOperator{\rank}{rank}
\DeclareMathOperator{\Pic}{Pic}
\newcommand\st{{\ \vert\ }}

\newtheorem{thm}{Theorem}

\title{Computing Fano 3-folds of index $\ge 3$}
\author{Gavin Brown
\and Kaori Suzuki}
\date{}

\makeatletter
\let\@old@maketitle=\@maketitle
\def\@maketitle{%
\@old@maketitle
{\def\thefootnote{}
\footnotetext{2000 \textit{Mathematics subject classification.} 
Primary 14J45 ; Secondary 14J30, 14Q15 \\ \hspace*{4mm} \textit{Key words and Phrases.} Fano 3-fold, Fano index, graded ring}}
}
\makeatother

\begin{document}
\maketitle

\begin{abstract}
We use the computer algebra system Magma to study graded rings
of Fano 3-folds of index $\ge 3$ in terms of their Hilbert series.
\end{abstract}

\section{Introduction}
Fano $3$-folds are, typically, the complex (projective) solution spaces
of homogeneous polynomial equations of low degree in 5 variables.
A quartic hypersurface is a classical example, for instance
\[
X_4 = (x_0^4 + \cdots + x_4^4=0) \subset \PP^4.
\]
In this example, the canonical class $K_{X_4}$ is represented simply
by a hyperplane section $A=(x_0=0)\subset X_4$, and so $X_4$ has
index (as defined below) equal to 1.
The cubic hypersurface $X_3 = (x_0^3 + \cdots + x_4^3=0) \subset \PP^4$
is also a Fano 3-fold, with $K_{X_3}=2A$ and so index~2.
Of course, there are more complicated examples involving more
variables, including weighted variables;  see \cite{IF} or \cite{ABR}
for an introduction to weighted projective space in this context.

By Suzuki \cite{Su}, the Fano index is bounded $f\le 19$
(and it does not take the values $12,14,15,16,18$).
We study Fano 3-folds of index $\ge 3$, especially the case $f=3$
generalising the conic hypersurface
$X_2 = (x_0^2 + \cdots + x_4^2=0) \subset \PP^4$;
see, for example, the lists of Iskovskikh and Prokhorov, \cite{IP}, Table~12.2.

Furthermore, in notation explained in the following section, we list
the number of possible numerical types (more precisely, of possible
Hilbert series) of Fano 3-folds of each index $f=3,\dots, 19$.
(The case of index $\ge 9$ is already proved in \cite{Su}.)
We work over the complex number field $\C$ throughout.

\begin{thm}
\label{thm!main}
For each $f=3,\dots, 19$, the the number of power series that could be
the Hilbert series of some $X,A$ with $X$ a Fano 3-fold of Fano
index~$f$ and $A$ a primitive ample divisor is:
\[
\begin{array}{c|cccccccccccc}
f &  3  & 4 &  5 &  6&   7  & 8  & 9  &10 & 11&  13&  17 & 19\\
\hline
\textrm{number of series}&231& 124 & 63 & 11 & 23&  10  & 2  & 1  & 3 & 2 &  1 &  1\\
\textrm{of which unstable}&50 &42 &29 & 5 &11 & 6 & 0 & 1 & 0 & 0 & 0 & 0
\end{array}
\]
The second line of this table indicates the number of these series that
cannot be realised by the Hilbert series of some Bogomolov--Kawamata stable 
Fano 3-fold $X,A$.
(See section~\ref{sec!alg}, Step 1(c)${}^+$, for a discussion of stability.
There are no Fano 3-folds of indices $f=12,14,15,16,18$.)
\end{thm}

Analogous methods for Fano 3-folds of index $\le 2$ work slightly
differently: in those cases there is another discrete invariant,
the genus, which does not play a role when $f\ge 3$.
This is why we stop here at $f=3$.
The following theorem is a result of our classification;
the proof is Step $2^+$ of section~\ref{sec!alg}.

\begin{thm}
\label{thm!KX}
$H^0(X,\Oh(-K_X))\not=0$ for any Fano 3-fold $X,A$ index $f\ge 3$.
\end{thm}

A first analysis of the possible realisations of these Hilbert series
in low codimension is in section~\ref{sec!egs} below.  As with all
the results in this paper, we used computer algebra---in our case,
the Magma system \cite{M}---in an essential way.  But this analysis,
and the list in codimension~4 especially, should be regarded only as
a list of possible examples and not a proved classification.
Tabulating these examples by codimension gives the following
(in which a blank entry is a zero);  all of these are stable.
\[
\begin{array}{c|cccccccccccc|c}
f & 3 & 4 & 5 & 6 & 7 & 8 & 9 & 10 & 11 & 13 & 17 & 19 & \textrm{total}\\
\hline
\textrm{codim $0$}& 0 & 1 & 1 & 0 & 1 & 0 & 0 & 0 & 1 & 1 & 1 & 1 & 7 \\
\textrm{codim $1$}& 7 & 2 & 5 & 1 & 4 & 3 & 2 & 0 & 2 & 1 &&& 27 \\
\textrm{codim $2$}& 6 & 7 & 1 & 0 & 0 &&&&&&&& 14 \\
\textrm{codim $3$}& 0 & 0 & 2 & 0 & 0 &&&&&&&& 2 \\
\textrm{codim $4$}& 3 & 2 & 1 & 1 & 1 &&&&&&&& 8
\end{array}
\]
Text files with the Magma code to make the classification of
Theorem~\ref{thm!main} and with all the proposed models is at the
webpage \cite{code}.

\section{Definitions and tools}

\paragraph{Basket of singularities}
Let the group $\Z/r$ of $r$th roots of unity act on $\C^3$ via the diagonal
representation $\eps\cdot(x,y,z)\mapsto (\eps^ax, \eps^ay, \eps^cy)$.
The (germ at the origin of the) quotient singularity $\C^3/(\Z/r)$
is denoted $\frac{1}r(a,b,c)$.  By Suzuki \cite{Su} Lemma~1.2,
when we work with Fano 3-folds of index~$f$ below, we may assume that
$b=-a$, $c=f$ and that $r$ is coprime to $a,b,c$.

We abbreviate the notation $\frac{1}r(a,-a,f)$ to $[r,a]$;  the index
$f$ is always clear from the context.
A {\em basket of singularities} is a collection (possibly with repeats)
of singularity germs $[r,a]$.

\paragraph{Fano 3-folds}
A {\em Fano 3-fold} is a normal projective 3-fold $X$ such that
(a) $-K_X$ is ample,
(b) $\rho(X) := \rank\Pic(X)= 1$, and
(c) $X$ has $\Q$-factorial
terminal singularities.
Without loss of generality, we may replace condition (c) by
the more restrictive condition: (c$^\prime$) $X$ is nonsingular apart
from a finite set of singularities equal to that of some basket.
(By Reid \cite{YPG} (10.2), this does not alter the Hilbert series
we compute and so our results hold as stated.  There may, however,
be series that are realised by Fano 3-folds satisfying (c) but not
(c$^\prime$)---but we do not know an example.)

The {\em Fano index $f=f(X)$} of a Fano 3-fold $X$ is
\[
f(X) = \max \{ m \in \Z_{> 0} \st -K_X = mA
\textrm{ for some Weil divisor }A \}
\]
where equality of divisors denotes linear equivalence of some multiple.
A Weil divisor $A$ for which $-K_X = fA$ is called
a {\em primitive ample divisor}.

\paragraph{Graded rings and Hilbert series}
A Fano $3$-fold $X$ with primitive ample divisor $A$, which we denote
by $X,A$ from now on, has a graded ring
\[
R(X, A) = \bigoplus_{n \geq 0} H^{0}(X, \Oh_{X}(nA)).
\] 
This graded ring is finitely generated, and
$X\cong\mathrm{Proj}~R(X,A)$.  The {\em Hilbert series $P_{X,A}(t)$
of $X,A$} is defined to be that of the graded ring $R(X,A)$:  thus
$\dim H^{0}(X, \Oh_{X}(nA))$ is the coefficient of $t^n$ in $P_{X,A}(t)$.

A choice of homogeneous generators for $R(X,A)$ determines a map
\begin{equation*}
X\hookrightarrow \PP^N = \PP(a_0,\dots,a_N)
\end{equation*}
into some weighted projective space (wps) $\PP^N$, where
$x_i\in H^{0}(X, \Oh_{X}(a_iA))$.
With this embedding for a minimal set of generators in mind,
we say that $X,A$ has {\em codimension} $N-3$.

\paragraph{The Riemann--Roch theorem}
Suzuki proves the appropriate version of Riemann--Roch in
this context, following Reid's plurigenus formula \cite{YPG},
to compute the dimensions of the graded pieces of $R(X,A)$.

For a singularity $p=\frac{1}{r}(a,-a,f)$ in $\cB$,
define $i_{p}(n) := -n/f\mod r$.  This always means
least residue modulo $r$, so that $0\le i_p(n) < r$.
When $r$ is clear from the context,
the notation $\overline{c}$ denotes the least residue of
$c$ modulo $r$.

\begin{thm}[\cite{Su} Theorem~1.4]
\label{thm:RR}
Let $X,A$ be a Fano 3-fold of index $f\ge 3$ and with basket $\cB$.
Then $p_n := \dim H^0(X,\Oh_{X}(nA))$ for any $n>-f$ is computed by
\begin{equation}
\label{eq!h0nA}
p_n = 1 + \frac{n(n+f)(2n+f)}{12}A^{3}
  + \frac{nAc_{2}(X)}{12} + \sum_{p=[r,a]\in\cB} c_p(n),
\end{equation}
where
$\displaystyle c_p(k)= -i_{p}(k)\frac{r^2-1}{12r} +
\sum_{j=1}^{i_{p}(k)-1}\frac{\overline{bj}(r-\overline{bj})}{2r}$
and $ab\equiv 1 \mod r$.

Summing these as a Hilbert series gives
\begin{eqnarray}
P_{X,A}(t) & = &
  \frac{1}{1-t} +
  \frac{(f^2+3f+2)t+(8-2f^2)t^2+(f^2-3f+2)t^3}{12(1-t)^4} A^{3} \nonumber \\
  && \qquad + \frac{t}{(1-t)^2}\frac{Ac_{2}(X)}{12} +
  \sum_{p\in\cB} \frac{1}{1-t^{r}}
    \sum_{k=1}^{r-1} c_p(k)t^k.			\label{eq!hs}
\end{eqnarray}
\end{thm}

Kawamata computes $Ac_2(X) = (1/f)(-K_Xc_2(X))$ in terms of $\cB$:
\begin{thm}[\cite{Ka}]
\label{thm:Ac2}
Let $X,A$ be a Fano 3-fold with basket $\cB$.  Then
\[
-K_Xc_2(X) = 24-\sum_{[r,a]\in\cB}\left(r-\frac{1}r\right).
\]
\end{thm}

\section{The algorithm for $3\le f\le 19$}
\label{sec!alg}

We explain our algorithm for arbitrary $3\le f\le 19$, and
we give explicit results only in the case $f=3$.

\paragraph{Step 1. Assembling possible baskets:}
A basket $\cB$ comprising germs $[r,a]$ of a Fano 3-fold must satisfy
several conditions.

\paragraph{Step 1(a)\quad Positive $Ac_2(X)$:}
Finiteness of the number is assured by Kawamata's condition
(\cite{Ka} in Theorem~2):
\[
-K_Xc_2(X)>0
\quad
\textrm{or equivalently}
\quad
\sum_{[r,a]\in\cB} \left(r - \frac{1}{r}\right) < 24.
\]
{\bf Result:}  2813 baskets satisfy Kawamata's condition.

\paragraph{Step 1(b)\quad Positive degree:}
The degree $A^3$ of $X,A$ can be computed from its basket $\cB$
by setting $n=-1$ in equation \eqref{eq!h0nA} since $H^0(X,\Oh(-A))=0$.
This degree must be strictly positive.

\noindent
{\bf Result:}  1295 of these baskets have $A^3>0$.

\paragraph{Step 1(b)$^+$\quad Excess vanishing:}
This condition can be strengthened since furthermore
$H^0(X,\Oh(nA))=0$ for each $n=-2,-3,\dots, -f+1$.
Enforcing this in equation \eqref{eq!h0nA}
has a significant effect once $f\ge 5$.

\paragraph{Step 1(c)\quad Bogomolov--Kawamata bound:}
By Suzuki \cite{Su} Proposition~2.4 and a consideration of the
stability of a tensor bundle in Kawamata \cite{Ka} Proposition~1,
\[
(4f^2-3f)A^3 \le 4f Ac_2(X).
\]
{\bf Result:}  231 of these baskets satisfy the Bogomolov--Kawamata bound.

\paragraph{Step 1(c)$^+$\quad Imposing stability:}
This is an optional step, and we do not include it in our full
classification.  It imposes the stronger condition
\[
f^2A^3 \le 3 Ac_2(X).
\]
Fano 3-folds (or their baskets) that satisfy this stronger bound are
called {\em Bogomolov--Kawamata stable},  being in the semistable part
of Kawamata's analysis \cite{Ka}.
While it is expected that this is the main case---possibly even
the only case---of the classification, this
condition is not known to hold for all Fano 3-folds.
All the examples we construct here are stable in this sense.

\noindent
{\bf Result:}  181 of these baskets are Bogomolov--Kawamata stable.

\paragraph{Step 2. Computing Hilbert series:}
For each basket in $\cB$, compute a power series $P(t)$ according
to the formula \eqref{eq!hs}.
By the expression of the formula, this is a rational function.
We also convert this into a power series (order 30 is sufficient for
our calculations); we use both representations later.

\paragraph{Step 2${}^+$. Sections of $-K_X$:}
Theorem~\ref{thm!KX} follows at once from the list of
Hilbert series.  We simply confirm that in each case the
coefficient of $t^f$ is nonzero.  Although we don't know
that each of these Hilbert series is realised by a Fano 3-fold,
certainly every Fano 3-fold (with $f\ge 3$) has Hilbert series
among our list.

\paragraph{Step 3. Estimating the degrees of generators:}
Suppose $P(t)=1+p_1t + p_2t^2 + \cdots$ is the Hilbert series of
some graded ring $R=\oplus_{d\ge 0} R_d$.  The following is a
standard method of guessing the degrees of some generators of a
minimal generating set of $R$.

Certainly $R$ must have $p_1$ generators of degree~1.
(Of course, this number may be zero.)
These generate at most a $q_2=\frac{1}2p_1(p_1-1)$-dimensional
subspace of $R_2$.
If $p_2-q_2\ge 0$, then $R$ must have at least $p_2-q_2$ generators
in degree~2.
On the other hand, if $p_2-q_2<0$, then this routine stops.
And so we continue into higher degree.

The calculation is made straightforward by the following observation.
If $n_1,\dots,n_d$ are the numbers of generators so far in degrees
$1$ up to $d$, then the number of monomials in degree $d+1$ they determine
(and so the maximum dimension space they could span in that degree)
is the coefficient of $t^{d+1}$ in the expansion
\[
\frac{1}{ (1-t)^{n_1} (1-t^2)^{n_2}\cdots (1-t^d)^{n_d} } = 1 + n_1t + \cdots.
\]
Such type changing (from rational functions to power series) is included
in most computer algebra systems, so this algorithm is easy to implement.

There are two important remarks.  First, the assumption of
generality (that the generators span a large space) can fail,
and this will change the degrees occurring in a minimal generating
set (although in small examples it will not reduce the number
of generators).  This is the main reason why our analysis is not
a complete proof, although it is compelling.

Second, in most cases this algorithm will not determine a complete
set of degrees for a minimal generating set.  This is the main
reason why we restrict our attention to low codimension when
proposing models, which we do next.

\paragraph{Step 4. Confirming small cases:}

The basket $\cB=\{ [2,1],[3,1],[7,3] \}$ with index $f=5$
determines the rational function
\[
P = \frac{t^8 + t^5 + t^4 + t^3 + 1}
{t^{13} - t^{12} - t^{11} + t^9 + t^8 - t^7 - t^6 + t^5 + t^4 - t^2 - t + 1}
\]
Expanded as a power series, this starts
\[
1 + t + 2t^2 + 4t^3 + 6t^4 + 9t^5 + 13t^6 + 18t^7 + 24t^8 + \cdots.
\]
The generator estimating routine above (called \verb!FindFirstGenerators(P)!
in Magma) predicts degrees $1,2,3,3,4,5$.
But $P$ is not of the form
\[
\frac{\textrm{polynomial in $t$}}
{\prod_{d=1,2,3,3,4,5} (1-t^d)}
\]
since the denominator still contains $t^6 + t^5 + t^4 + t^3 + t^2 + t + 1$.
The solution is clear:  include $7$ as the degree of a generator.
From the Hilbert series point of view, this absorbs the excess
factor in the denominator;  from the basket point of view, this
provides the cyclic group action to generate the contribution of
the quotient singularity $[7,3]=\frac{1}7(3,8,5)$ in the basket.

The final form of the Hilbert series is thus
\[
\frac{-t^{20} + t^{14} + t^{13} + t^{12} + t^{11} - t^9 - t^8 - t^7 - t^6 + 1}
{\prod_{d=1,2,3,3,4,5,7} (1-t^d)}
\]
which suggests a variety defined by 5 equations
of weights $6,7,8,9,10$:
\[
X_{6,7,8,9,10} \subset \PP^6(1,2,3,3,4,5,7).
\]
In fact, these equations can be written as the five maximal Pfaffians of
a skew $5\times 5$ matrix, as in \cite{ABR} Remark~1.8 or \cite{kinosaki}
section~4, and it can be checked that this $X$ is a Fano 3-fold with
singularities equal to the basket.

\section{Classification in low codimension}
\label{sec!egs}

We distinguish between cases in codimension $\le 3$, where
we can write down equations of Fano 3-folds and check
their properties explicitly, and codimension~4, where
calculations are more difficult.
Tables of these results are given below,
and the webpage \cite{code} contains these and all other
Hilbert series as Magma output, as well as the Magma code to
generate them.

\paragraph{Examples in codimension at most 3}

Only seven weighted projective spaces are themselves are Fano 3-folds.
These are:
$\PP^3$ with $f=4$;
$\PP(1,1,1,2)$ with $f=5$;
$\PP(1,1,2,3)$ with $f=7$;
$\PP(1,2,3,5)$ with $f=11$;
$\PP(1,3,4,5)$ with $f=13$;
$\PP(2,3,5,7)$ with $f=17$;
$\PP(3,4,5,7)$ with $f=19$.

For hypersurfaces or in codimension~2, listed in Tables~\ref{tab!cod1}
and \ref{tab!cod2},
the equations are simply generic polynomials of the indicated
degrees.  Table~\ref{tab!cod3} lists those in codimension~3;
here one must build a $5\times 5$ skew
matrix of forms (as in \cite{ABR} Remark~1.8),
and then the equations are its five $4\times 4$ Pfaffians.
It is a mystery why there are so few families here for $f\ge 3$;
by comparison, in the case $f=1$ there are 70 families in
codimension~3.

\paragraph{Examples in codimension 4 are more subtle}
The Hilbert series routines and guesses of additional weights
work in exactly the same way in codimension~4 as in lower
codimension.  But it is not easy to write down an example
of a ring with given generator degrees in codimension~4.  In other
graded ring calculations, such as for K3 surfaces in \cite{B},
there is much use of projection and unprojection methods.
But (Gorenstein) projection of a Fano of higher index
does not result in another Fano.  Nevertheless, the projection
construction of a K3 surface section $S=(x=0)\subset X$,
where $x$ is a variable in degree~$f$, can be a guide.
We propose the list of examples in Table~\ref{tab!cod4},
although none has been constructed explicitly.
As justification, we give an example to illustrate what goes
wrong with the possible codimension~4 models that we have
rejected---the proposals listed in Table~\ref{tab!cod4} are
exactly those candidates that do not suffer from this obstruction.

Let index $f=4$ and basket $\cB=\{ [5,2] \}$;  these (stable) data
determine a Hilbert series $P(t)$.  Suppose we can construct a
Fano 3-fold $X,A$ having Hilbert series $P$.
Considerations as above suggest the degrees of a minimal
set of eight generators for the ring $R(X,A)$ could be
$1,1,1,2,2,3,4,5$ so that $X$ is in codimension~4.
And indeed there is a family of codimension~4 K3 surfaces in
$\PP^6(1,1,1,2,2,3,5)$ that could be the K3 sections
$(t=0)\subset X$, where $t$ is the variable on $X$ of weight~4.
Now a typical such K3 surface $S$ admits a projection to a K3
surface of codimension~3 in $\PP^5(1,1,1,2,2,3)$---this is
simply the elimination of the degree~5 variable from the
ideal defining $S$ (using the Groebner basis with respect to a
standard lexicographic monomial order with $t$ big, for instance).
The image is in codimension~3, and its equations are the five
Pfaffians of a skew $5\times 5$ matrix of forms.  Crucially,
one calculates that the forms appearing here
each have degree $\le 3$.
So the analogous projection of $X$ would have
equations that not involving the variable $t$, and this
would force a non-terminal singularity onto $X$ itself.

\noindent
Gavin Brown\\
IMSAS, University of Kent, CT2 7AF, UK.  Email: \texttt{gdb@kent.ac.uk}

\vspace{3mm}
\noindent
Kaori Suzuki\\
Tokyo Institute of Technology, 2-12-1 Ookayama, Meguro-ku, 152-8550, Japan.
Email:~\texttt{k-suzuki@math.titech.ac.jp}

\pagebreak

\begin{table}[h]
\[
\begin{array}{c|c|c|c|c}
f & \textrm{Fano 3-fold $X\subset\PP^4$} & A^3 & Ac_2(X) & \textrm{Basket $\cB$}  \\
\hline
\hline
3 & X_2\subset\PP^4		& 2 & 8 & \textrm{no singularities}	\\
  & X_3\subset\PP(1,1,1,1,2)	& 3/2 & 15/2 & [2,1]	\\
  & X_4\subset\PP(1,1,1,2,2)	& 1 & 7/12 & 2\times[2,1]	\\
  & X_6\subset\PP(1,1,2,2,3)	& 1/2 & 13/2 & 3\times[2,1]	\\
  & X_{12}\subset\PP(1,2,3,4,5)	& 1/10 & 49/10 & 3\times[2,1],[5,1]	\\
  & X_{15}\subset\PP(1,2,3,5,7)	& 1/14 & 73/14 & [2,1],[7,2]	\\
  & X_{21}\subset\PP(1,3,5,7,8)	& 1/40 & 151/40 & [5,2],[8,1]	\\
\hline
4 & X_4\subset\PP(1,1,1,2,3)	& 2/3 & 16/3 & [3,1]	\\
 & X_6\subset\PP(1,1,2,3,3)	& 1/3 & 14/3 & 2\times[3,1]	\\
\hline
5 & X_4\subset\PP(1,1,2,2,3)	& 1/3 & 11/3 & 2\times[2,1],[3,1]	\\
 & X_6\subset\PP(1,1,2,3,4)	& 1/4 & 15/4 & [2,1],[4,1]	\\
 & X_6\subset\PP(1,2,2,3,3)	& 1/6 & 17/6 & [2,1]	\\
 & X_{10}\subset\PP(1,2,3,4,5)	& 1/12 & 35/12 & 2\times[2,1],[3,1],[4,1]	\\
 & X_{15}\subset\PP(1,3,4,5,7) & 1/28 & 75/28 & [4,1],[7,3] \\
\hline
6 & X_6\subset\PP(1,1,2,3,5) & 1/5 & 16/5 & [5,2]	\\
\hline
7 & X_6\subset\PP(1,2,2,3,5) & 1/10 & 21/10 & 3\times[2,1],[5,2]	\\
  & X_6\subset\PP(1,2,3,3,4) & 1/12 & 23/12 & [2,1],2\times[3,1],[4,1]	\\
  & X_8\subset\PP(1,2,3,4,5) & 1/15 & 29/15 & 2\times[2,1],[3,1],[5,1]	\\
  & X_{14}\subset\PP(2,3,4,5,7) & 1/60 & 71/60 & 3\times[2,1],[3,1],[4,1],[5,2]	\\
\hline
8 & X_{6}\subset\PP(1,2,3,3,5) & 1/15 & 26/15 & 2\times[3,1],[5,2]	\\
  & X_{10}\subset\PP(1,2,3,5,7) & 1/21 & 38/21 & [3,1],[7,2]	\\
  & X_{12}\subset\PP(1,3,4,5,7) & 1/35 & 54/35 & [5,1],[7,3]	\\
\hline
9 & X_6 \subset \PP(1,2,3,4,5) & 1/20 & 31/20 & [2,1],[4,1],[5,2]	\\
  & X_{12} \subset \PP(2,3,4,5,7) & 1/70 & 61/70 & 3\times[2,1],[5,2],[7,3]	\\
\hline
11 & X_{12} \subset \PP(1,4,5,6,7) & 1/70 & 69/70 & [2,1],[5,1],[7,1]	\\
   & X_{10} \subset \PP(2,3,4,5,7) & 1/84 & 59/84 & 2\times[2,1],[3,1],[4,1],[7,2]	\\
\hline
13 & X_{12} \subset \PP(3,4,5,6,7) & 1/210 & 89/210 & [2,1],2\times[3,1],[5,1],[7,3]	\\
\end{array}
\]
\caption{Fano 3-folds in codimension 1\label{tab!cod1}}
\end{table}

\begin{table}[h]
\[
\begin{array}{c|c|c|c|c}
f & \textrm{Fano 3-fold $X\subset\PP^5$} & A^3 & Ac_2(X) & \textrm{Basket $\cB$} \\
\hline
\hline
3 & X_{6,6}\subset\PP(1,1,2,3,3,5) & 2/5 & 32/5 & [5,2]
\\
 & X_{6,6}\subset\PP(1,2,2,3,3,4) & 1/4 & 19/4 & 4\times [2,1], [4,1]
\\
 & X_{6,9}\subset\PP(1,2,3,3,4,5) & 3/20 & 93/20 & [2,1], [4,1], [5,2]
\\
 & X_{12,15}\subset\PP(1,3,4,5,6,11) & 1/22 & 85/22 & [2,1], [11,5]
\\
 & X_{9,12}\subset\PP(2,3,3,4,5,7) & 3/70 & 183/70 & 3\times [2,1], [5,2], [7,3]
\\
 & X_{12,15}\subset\PP(3,3,4,5,7,8) & 1/56 & 103/56 & [4,1], [7,3], [8,3]
\\
\hline
4 & X_{6,8}\subset\PP(1,2,3,3,4,5) & 2/15 & 52/15 & 2\times[3,1],[5,2]
\\
 & X_{8,10}\subset\PP(1,2,3,4,5,7) & 2/21 & 76/21 & [3,1],[7,2]
\\
 & X_{8,12}\subset\PP(1,3,4,4,5,7) & 2/35 & 108/35 & [5,1],[7,3]
\\
 & X_{10,12}\subset\PP(1,3,4,5,6,7) & 1/21 & 62/21 & 2\times[3,1],[7,1]
\\
 & X_{10,12}\subset\PP(2,3,4,5,5,7) & 1/35 & 66/35 & 2\times[5,2],[7,2]
\\
 & X_{12,14}\subset\PP(2,3,4,5,7,9) & 1/45 & 86/45 & [3,1],[5,2],[9,2]
\\
 & X_{18,20}\subset\PP(4,5,6,7,9,11) & 1/231 & 206/231 & [3,1],[7,2],[11,5]
\\
\hline
5 & X_{10,15}\subset\PP(2,3,5,5,7,8) & 1/56 & 87/56 & [2,1],[7,2],[8,3]
\end{array}
\]
\caption{Fano 3-folds in codimension 2\label{tab!cod2}}
\end{table}

\begin{table}
\[
\begin{array}{c|c|c|c|c}
f & \textrm{Fano 3-fold $X\subset\PP^6$} & A^3 & Ac_2(X) & \textrm{Basket $\cB$} \\
\hline
\hline
5 & X_{6,7,8,9,10}\subset\PP(1,2,3,3,4,5,7) & 5/42 & 109/42 & [2,1], [3,1], [7,3]
\\
 & X_{12,13,14,15,16}\subset\PP(1,4,5,6,7,8,9) & 1/36 & 71/36 & [2,1], [4,1], [9,1]
\end{array}
\]
\caption{Fano 3-folds in codimension 3\label{tab!cod3}}
\end{table}

\begin{table}
\[
\begin{array}{c|c|c|c|c}
f & \textrm{Fano 3-fold $X\subset\PP^7$} & A^3 & Ac_2(X) & \textrm{Basket $\cB$} \\
\hline
\hline
3 & X\subset\PP(1,1,1,1,1,3,3,4) & 9/4 & 27/4 & [4,1]
\\
 & X\subset\PP(1,1,2,2,3,3,4,5) & 3/5 & 27/5 & [2,1], [2,1], [5,1]
\\
 & X\subset\PP(1,1,2,2,3,3,5,7) & 4/7 & 40/7 & [7,2]
\\
4 & X\subset\PP(1,1,2,2,3,3,4,5) & 7/15 & 62/15 & [3,1], [5,2]
\\
 & X\subset\PP(1,1,2,2,3,4,5,7) & 3/7 & 30/7 & [7,2]
\\
5 & X\subset\PP(1,1,2,3,3,4,5,7) & 2/7 & 24/7 & [7,3]
\\
6 & X\subset\PP(1,2,3,4,5,5,6,7) & 3/35 & 72/35 & [5,2], [7,2]
\end{array}
\]
\caption{Proposals for Fano 3-folds in codimension 4\label{tab!cod4}}
\end{table}

\end{document}